\documentclass{article}
\usepackage[utf8]{inputenc}
\usepackage[english]{babel}
\usepackage{amsmath,amssymb}
\usepackage{fancyhdr}
\usepackage{authblk}

\usepackage[a4paper]{geometry}

\usepackage{eurosym}
\usepackage{algorithm}
\usepackage{algpseudocode}

\usepackage{comment}

\usepackage{graphicx}
\usepackage{subcaption}

\usepackage{xcolor}
\usepackage[normalem]{ulem}

\usepackage{natbib}
\setcitestyle{authoryear,open={(},close={)}} 
\usepackage{url}

\newtheorem{theorem}{Theorem}

\newtheorem{remark}[theorem]{Remark}

\newcommand{\ECTRL}{\textsc{Eurocontrol}}

\title{Data-driven optimization for Air Traffic Flow Management with trajectory preferences}
\author[a]{Luigi De Giovanni}
\author[b,c]{Guglielmo Lulli}
\author[c]{Carlo Lancia}

\affil[a]{Dipartimento di Matematica “Tullio Levi-Civita”, Università degli Studi di Padova, Padova, Italy -- luigi@math.unipd.it}
\affil[a]{Management Science, Lancaster University, Lancaster, United Kingdom,  -- g.lulli@lancaster.ac.uk}
\affil[c]{Dip. di Informatica, Sistemistica e Comunicazione, Universit\`a di Milano-Bicocca}

\date{}

\begin{document}

\maketitle

\begin{abstract}
In this paper, we present a novel data-driven optimization approach for trajectory based air traffic flow management (ATFM). 
A key aspect of the proposed approach is the inclusion of airspace users' trajectory preferences, which are computed from traffic data by combining clustering and classification techniques. Machine learning is also used to extract consistent trajectory options, while optimization is applied to resolve demand-capacity imbalances by means of a mathematical programming model that judiciously assigns a feasible 4D trajectory and a possible ground delay to each flight. 
The methodology has been tested on instances extracted from real air traffic data repositories. With more than 32,000 flights considered, we solve the largest instances of the ATFM problem available in the literature in short computational times that are reasonable from the practical point of view. As a by-product, we highlight the trade-off between preferences and delays as well as the potential benefits. Indeed, computing efficient solutions of the problem facilitates a consensus between network manager and airspace users. In view of the level of accuracy of the solutions and of the excellent computational performance, we are optimistic that the proposed approach may provide a significant contribution to the development of the next generation of air traffic flow management tools.
\end{abstract}

\textbf{keywords:} Air Traffic Flow Management, Optimization, Machine Learning, Preferences.



\section{Introduction}
\label{sec::intro}

The air traffic industry is a strategically important sector that makes crucial contribution to the overall world economy and employment. Before the COVID-19 epidemic,  in Europe, aviation was supporting almost 5 million jobs and contributed \euro 300B to the European GDP (EU Commission). Despite its vital role, the air traffic system was experiencing severe congestion phenomena on a daily base that were impairing the air traffic industry's sustainability. \ECTRL~ reported 19.2 mln minutes of en-route air traffic flow management delays for year 2018. The total air traffic delays -- in the same calendar year -- costed the EU economy \euro 17.6B ~\citep{PRR2018}. This situation was even more exacerbated in the US with overall delay costs for year 2018 estimated in \$28B. 

Although  air travel is currently one of the economy's sectors most affected by the COVID-19 epidemic \citep[see for instance the report of][]{sp2020}, there is a general agreement among experts that the sector will recover. Indeed, all the major R\&D programmes for modernizing the air transport systems around the world are continuing with their plan and long-term objectives. One of the cornerstones of these initiatives is the implementation of the ICAO Trajectory Based Operations (TBO) concept. The ambition of the TBO concept is to provide the capability of flying a path that is as close as possible to the user-preferred one. From a full implementation of the TBO concept, airspace users expect to achieve a higher degree of flexibility to manage their operations and meet their business objectives. To enable TBO operations, it is critical to resolve demand-capacity imbalances in a more efficient way than today practise. 

Several air traffic flow management (ATFM) initiatives and dedicated tools are currently used to resolve air traffic  congestion. For instance, in Europe, the Network Manager -- i.e., the authority in charge of smoothing traffic flows on the entire air traffic network -- allocates time slots by the Computer-Assisted Slot Allocation (CASA) algorithm \citep[e.g., see][]{casa}. CASA implements a set of rules with the First Planned First Served (FPFS) as the primary one.  Indeed, the FPFS rule (also known as ration-by-schedule mechanism) underpins all the major ATFM tools and initiatives in the US as well.  However, the CASA algorithm, as well as most if not all of the ATFM procedures deployed so far, is far from being optimal, with the assignment of ``unnecessary'' large delays as highlighted in~\cite{estes2019} and~\cite{ruiz2019}. This undesirable behaviour is due  to the application of the FPFS rule in the decision process, without any estimate of the potential downstream effects of such decisions. 

To address the optimality issue, meaning of computing  (control) decisions that optimize the system performance, the scientific community has developed several mathematical models and optimization algorithms. Nevertheless, most of these approaches do not provide representation of flight trajectories at the level of accuracy that is needed for TBO operations. TBO requires operational trajectories in the four dimensions, i.e., space -- including altitude -- and time. Moreover, there are still some intrinsic issues on the computational side as well as on their impact to practise. On the computational side, although remarkable progresses have been made on developing faster algorithms, the scalability problem is still a concern. Indeed, real daily instances of the problem with 4D trajectories have not been solved yet. 

As far as the impact to practise, it is still quite limited because the proposed solutions often fail to meet airspace users' business objective and operational requirements. The importance of fulfilling the business objectives and operational requirements is not only desirable, it is now not deferable. For the sake of the truth, this failure is not only due to modelling capabilities but for the most part to a lack of information. To enable a collaborative decision-making environment in air traffic flow management with a more active engagement of airspace users in the decision making process, the Collaborative Trajectory Options Program (CTOP) has been developed and deployed in the US \citep{ctop14}. For each flight, airspace users submit a set of feasible trajectories each with a relative cost. The relative trajectory cost is an expression of preference and relative value. If congestion occurs, the Federal Aviation Administration (FAA) assigns to each flight the most desirable trajectory among the feasible ones, following the FPFS rule. However, quoting \cite{hoff2018}, ``\textit{despite CTOP having been deployed in 2014, the FAA has been reluctant to use CTOP because almost none of the air carriers are prepared to generate} [and manage] \textit{trajectory option sets.}'' The lack of tools for generating and managing trajectory option sets combined to the lack of incentives for disclosing preference information \citep[][]{estes2019}  proved to be a barrier to the adoption of CTOP by airspace users.

The concept of trajectory preference has been also widely discussed by European practitioners \cite[see, e.g.,][]{workshop}. A preference was conceived as a partial order of feasible control options -- i.e., time, flight level, lateral deviation or a combination of them --  to deviate a flight from the preferred  trajectory in order to manage delays at tactical level. These options are all observable features of a trajectory that determine a preference. However, determinants of a preference are in many cases only partially known or, even, unknown. 

The ability of capturing airspace users' preferences combined with the selection and assignment of feasible 4D trajectories has the potential of  facilitating a consensus of all the parties involved on the identified ATFM decisions thus accelerating the ATFM decision process. In this respect, the preferences' criterion has to be combined with the criterion of delays, traditionally considered to assess the performance of ATFM initiatives. We here suggest  computing the Pareto efficient solution that minimize the total delay, because the proposed approach  is intended for the ATFM function and the ATFM authority more in general, e.g., \ECTRL's Network Manager.

To overcome the lack of information on trajectories' preference, we use machine learning tools. More in particular, we first extract consistent flight plans (trajectories) from historical data, and cluster them in terms of geometry and speed. We then use a random forest to predict the cluster membership of flown trajectories given relevant features of the flight, which provides, as a byproduct, a measure of the trajectory (group) preference for the flight itself. This information feeds a mathematical programming model that assigns a trajectory and a possible ground delay to each flight, and optimally resolves demand-capacity imbalances. The proposed approach has been tested on instances extracted from the \ECTRL~ data repository, and provides a large real example of Artifical Intelligence \& Optimization in the air traffic domain. Indeed, with more than 32,000 flights considered, we solved the largest instances of the ATFM problem available in the literature. These instances also enjoy the distinguishing feature that both the air traffic demand and the air traffic system, i.e., the network of en-route sectors and airports and corresponding capacities, are based on real data, which is unique in the ATFM literature at least at this level of size. To increase the level of fidelity of the proposed solutions, we accurately model airspace capacity restrictions as implemented in practise.

In summary, the main contributions of the paper are: 
{\it i)} the use of machine learning techniques to compute trajectory options and related airspace users' preferences; 
{\it ii)} the development of a mathematical model to facilitate the ATFM planning process that incorporates preferences and a realistic representation of capacity constraints;
{\it iii)} the design of effective algorithms to solve the mathematical model for real size instances; {\it iv)} the experimentation of the proposed approach on large real daily instances of the problem, highlighting the computational challenges as well as the benefits of the proposed approach.

The remainder of the paper is organized as follows. Section~\ref{sec::literature} reviews the relevant literature on the problem herein addressed. Section~\ref{sec::method} presents the details of the proposed methodological approach, including the machine learning models for trajectories and trajectories' preference extrapolation, and the mathematical formulation of the trajectory based ATFM. The  computational approach to solve the proposed model for large size instances of the problem is discussed in Section \ref{sec::soltech}. Section \ref{sec::CompExp} is devoted to computational experiments. It describes in detail the realistic daily instances extracted from \ECTRL~ data repository, the results and the analysis of solutions. Finally, Section~\ref{sec::concl} concludes the paper.

\section{Literature review}
\label{sec::literature}

In the last four decades, the ATFM domain has continuously attracted the interest of the research community with the aim of finding effective and efficient ways to resolve unbalances between air traffic demand and air system capacity. The early focus of research activities was on airports' congestion. This effort led to a better understanding and improvement of Ground Delay Programs as an effective tool to control and manage the temporary excess of air traffic demand at airports, see~\cite{vossen2011} for a detailed survey on the topic. In more recent years, the scope of ATFM initiatives, and consequently of mathematical models and algorithms, enlarged to resolve en-route congestion as well. For this class of models, as highlighted by \cite{LOd_TS07}, the optimal ATFM strategies can be complex and occasionally counterintuitive due to the simultaneous presence of airports and en-route sectors capacity constraints. This explains the clear  need of holistic approaches to the ATFM problem, which is indeed ubiquitous in the most recent ATFM models and/or formulations. The class of models that are relevant for the research and the methodological approach proposed in this paper are the so-called  Lagrangian models, i.e., models that reproduce the full trajectory of each flight.  

At the best of our knowledge, one of the first attempts to consider both airports and en-route sectors' capacities in a deterministic trajectory-based setting was proposed by \cite{lindsay93}. The authors formulated a 0-1 integer programming model for assigning ground and airborne holding delays to single flights in the presence of both airport and airspace capacity constraints. \cite{BerS98} provided a formulation of a similar model. The main feature of the proposed formulation is its tightness. Indeed, several of the constraints define facets of the polyhedron of solutions.  As a result of the good mathematical structure, the Bertsimas-Stock formulation enables the fast computation of optimal solutions. In a second paper, \cite{BerS200} developed a dynamic, multicommodity, integer network-flow model to include rerouting as a control option. However, the proposed approach did not  prove to be very effective in solving large instances of the problem. Issue that was addressed by \cite{BLO}. This model includes most of the ATFM control options while preserving the good mathematical structure of \cite{BerS98}. The main innovative feature of the model is the formulation of ``rerouting decisions'' in a very compact way. The authors were able to solve large instances of the problem with more than 6,500 flights. \cite{augustin2012} developed a model for ATFM using the same concept of rerouting as in \cite{BLO}. 

A different modeling approach is the ``trajectory" formulation first proposed in the ATFM context by \cite{richard2011}. Given the large number of variables of these formulations, they are usually solved by means of a column generation approach. \cite{hamsa2014} designed an algorithm for this class of formulations that is computationally very efficient with good scalability properties. The algorithm solved NAS-wide numerical examples with up to 17,500 flights in short computational times.
Somehow of a similar nature, it is the model proposed by \cite{sherali2002}, which also assigns flights to trajectories. In this case, the set of feasible 4D trajectories is explicitly given as input of the model. To the best of our knowledge, this was the first model to consider  trajectories in 4D, which is crucial for TBO operations. Indeed, all other models listed above describe the flights' trajectories in a two dimensional geographical space, meaning that no flight level information is attached to the trajectories. More recently,  \cite{dalsasso2018_TRB, dalSasso_ejor} explicitly modeled flight level information in a compact formulation. 
The authors highlighted some of the challenges in modeling 4D trajectories in terms of producing acceptable solutions for airspace users. These two papers do also represent the first attempt to explicitly include preferences in a realistic setting of the ATFM domain. To overcome the issue of information sharing, the authors developed a multi-objective binary integer programming model with the ambition of  highlighting the trade offs involved with the identified primitives of  airspace users' preferences, i.e., delay, flight level and route charges.

\section{Methodology}
\label{sec::method}

The key aspect of the methodology developed in this paper is the combination of machine learning and optimization techniques to compute viable ATFM solutions for real instances of the problem at a continental scale. In particular, machine learning methods are used to mine trajectory options  (\S \ref{sec::trjclas}) and related airspace users preferences (\S \ref{sec::pref}). An optimization model is then used to assign a feasible trajectory to each flight and resolve demand-capacity imbalances (\S \ref{sec::formulation}).

\subsection{Trajectory extrapolation and classification}
\label{sec::trjclas}

Our methodology stems from the availability of consolidated historical data on filed flight plans stored in data repositories such as \ECTRL~Demand Data Repository (DDR2). In this work, we consider the tactical flight plan (TFP), i.e., the last filed flight plan before the application of any operational regulation. The TFP better represents the airspace user's choice and preference at tactical level, which is relevant for the ATFM problem herein addressed.  
Each flight plan (trajectory) is described as a sequence of either waypoints (fixes) or sector entering points. For each data point in the sequence, the following information is available: latitude, longitude, flight level and  time-of-the-day.
Given a pair $(i,j)$ of origin and destination airports, we extract a consistent set of trajectories using Algorithm~\ref{alg:extrTraj}.

\begin{algorithm}
\caption{\textsc{Trajectory Extrapolation} }
\label{alg:extrTraj}
\footnotesize

{\bf Input:} origin airport $i$, destination airport $j$, filed TFPs.

{\bf Output:} set of trajectory options between $i$ and $j$ and related clustering.
\begin{algorithmic}[1]
\State \label{a:te:extract}
\textsc{Extract filed TFPs} between $i$ and $j$ in terms of waypoints.
\State \label{a:te:resample}
\textsc{Re-sample trajectories} from differently many waypoints to a same number of equidistant 4D points.
\State \label{a:te:clean}
\textsc{Trim trajectories} to remove takeoff and approach sample points.
\State \label{a:te:cluster}
\textsc{Apply a clustering procedure} to compute clusters of trajectories and detect outliers.
\State \label{a:te:remove}
\textsc{Remove outliers} to obtain clustered trajectories between $i$ and $j$. 
\end{algorithmic}
\end{algorithm}

The core of the procedure is  a density-based clustering of filed trajectories. The method is inspired by the one proposed by \cite{gariel2011} and adopted by, e.g., \cite{liu2021}. Density-based clustering groups similar trajectories according to  their pair-wise distance in a given metric space. It has the advantage of not requiring an \emph{a priori} knowledge of the number of clusters and it is robust to outliers, i.e., peculiar trajectories that are markedly different from any other trajectory in the dataset. Outliers' detection is relevant to our choice model, because they correspond to trajectories determined by factors that are hardly related to preferences (unusual weather conditions, strikes, adverse natural events etc.).

Algorithm \ref{alg:extrTraj} starts by extracting all the filed TFPs related to flights between airports $i$ and $j$ operated in a relevant time interval, e.g., a season.  
Trajectories are extracted in terms of waypoints and they are made of sequences of different length. Before the clustering step, we re-sample trajectories from differently many waypoints to same number $n$ of 4D points, where $n$ is twice the length of the longest sequence extracted between $i$ and $j$ (Step \ref{a:te:resample}). Each sample point consists of longitude, latitude, flight level and elapsed time from departure: latitude and longitude are computed in such a way that consecutive sample points are equidistant in the latitude-longitude plane, while flight level and elapsed time are obtained by linear interpolation. After re-sampling, every trajectory 
is described by a vector in $\mathbb{R}^{4n}$.
Min-max scaling is used to map each of the $4n$ features of a trajectory $x$ onto the interval [0, 1]. 
In order to remove segments related to takeoff and approach manoeuvres, that are less relevant at the pre-tactical stage, we trim off trajectory head and tail, approximated by the initial and final 10\% points. 
Step \ref{a:te:cluster}  
clusters trajectories into homogeneous groups using DBSCAN with the cosine distance metric, which is defined as
\[
  1 - \frac{\sum_{i=1}^{4n} x_i\,y_i}{\sqrt{\sum_{i=1}^{4n} x_i^2}\sqrt{\sum_{i=1}^{4n} y_i^2}}\,,
\]
for any two trajectories $x, y \in \mathbb{R}^{4n}$. 
Finally, we discard outliers to determine the set of consistent trajectory options for flights between airport $i$ and airport $j$.

We remark that the proposed methodology not only provides a set of trajectories that can be consistently flown between airport pairs, but also bins trajectories into clusters based on their geometry (including the flight level) and operating speed. Since trajectories in the same cluster are similar to each other, we assume that they share the same preference scores, and the obtained clusters provide the basis for preference extrapolation.

\subsection{Flight classification and preference extrapolation}
\label{sec::pref}

Towards preference-aware ATFM, we need a measure of the preference that an airspace user assigns to a trajectory extracted by Algorithm \ref{alg:extrTraj}.
We assume that this preference depends on the flight features and, according to the closing remark of \S\ref{sec::trjclas}, on the cluster of the trajectory.
Learning flight preferences can be thus reduced to learning how the features of a flight are related to the cluster of the trajectory it flies.

To this end, we develop a machine learning tool, which is able to predict the (preferred) cluster based on the flight features. This is the first step of Algorithm \ref{alg:fliPref}, that we propose for preference extrapolation.
\begin{algorithm}
\caption{\textsc{Computation of Flights' Preferences} }
\label{alg:fliPref}
\footnotesize
{\bf Input:} origin and destination airports, flights with features and cluster of the filed TFP.

{\bf Output:} preference of each flight for each trajectory option.
\begin{algorithmic}[1]
\State 
\label{a:pe:learn-tree}\textsc{Train a Random Forest Classifier} that predicts the cluster of the filed TFP, based on the flight features.
\For{each flight}
    \State \label{a:pe:classify}
    \textsc{Classify} the flight and let each tree in the Random Forest estimate the associations to all of the clusters.
    \State  \label{a:pe:pref-function}
    \textsc{Compute preferences} of the flight by averaging individual-tree associations. 
\EndFor
\end{algorithmic}
\end{algorithm}
In particular, we train a random forest classifier on the following flight data: day of the week, week number (for seasonal effects), anonymized airline code, airline type (legacy/low cost), and aircraft model. All variables are considered as categorical and binarized by one-hot encoding. We recall that a random forest classifier is an ensemble of binary trees. 

The internal nodes of each tree represent a boolean condition on a predictive variable, while its leaves measure the \emph{association} to each cluster as the proportion of samples of that cluster falling into the same leaf. In a random forest, each tree is fit on a subset of the predictive features and a subset of the training instances (sampled with replacement). Each constituent of the ensemble classifies a flight, determine associations and votes for one of the clusters previously extracted (step \ref{a:pe:classify}) based on maximum association. In our study, associations are also used to estimate the preference of a flight for trajectories in each cluster as the average associations from the trees in the ensemble (step \ref{a:pe:pref-function}).


Table~\ref{tab:results-rf} shows the performance of the classifier on predicting the cluster of the flown trajectory from the set learned via Algorithm~\ref{alg:fliPref} for some origin-destination pairs. In addition to the number of trajectories, outliers and clusters, we report the weighted average of f1-score, precision, and recall, computed using a nested cross-validation scheme with a 10-fold stratified split \citep[see, e.g.,][]{zheng2015}. The average of each metric is computed over 10 runs, and we also report the resulting standard error of the mean (\emph{sem}).
\begin{table}
\begin{small}
\centering
\caption{Cross-validated prediction performance of the random forest classifier.}
\begin{tabular}{lccccc} \hline \hline
Origin-Destination & \#traj (\#outl) & \#clusters & F1-score (sem) & Precision (sem) & Recall (sem) \\ \hline 
Amsterdam - London &  4184 (131) & 4 & 0.87 (0.007) & 0.88 (0.006) & 0.86 (0.008)\\
Frankfurt - Athens &  338 (22) & 3 & 0.96 (0.010) & 0.96 (0.011) & 0.96 (0.009)\\
Istanbul - Frankfurt &  1112 (159) & 5 & 0.89 (0.009) & 0.90 (0.010) & 0.89 (0.008)\\
London - Athens &  979 (51) & 6 & 0.90 (0.007) & 0.90 (0.007) & 0.89 (0.008)\\
London - Frankfurt &  2096 (142) & 8 & 0.86 (0.007) & 0.87 (0.007) & 0.86 (0.008)\\
London - Istanbul &  1431 (102) & 4 & 0.76 (0.013) & 0.77 (0.013) & 0.76 (0.013)\\
Madrid - London &  2418 (157) & 4 & 0.95 (0.004) & 0.95 (0.004) & 0.95 (0.004)\\
Rome - Barcelona &  1071 (69) & 3 & 0.82 (0.014) & 0.82 (0.015) & 0.82 (0.013)\\
Rome - London &  2168 (161) & 6 & 0.80 (0.006) & 0.81 (0.006) & 0.81 (0.006)\\
Rome - Paris &  1810 (155) & 7 & 0.96 (0.004) & 0.96 (0.004) & 0.97 (0.003)\\
\hline \hline
    \multicolumn{6}{l}{\footnotesize sem: standard error of the mean}
    \end{tabular}
    \label{tab:results-rf}
\end{small}
\end{table}
The performance of the classifier are commonly good, showing that the approach is capable of mapping to a satisfactory extent flight characteristics to airlines preferences.

\subsection{An Integer Linear Programming Model for trajectory selection} 
\label{sec::formulation}

Once  the set of eligible 4D trajectories is retrieved, the goal of our approach is to solve the ATFM problem by a mathematical model that assigns one trajectory to each flight, together with a possible ground delay to meet capacity restrictions. Other tactical control options -- e.g., flight level capping -- are directly embedded in the trajectories, whereas we do not model control options of a more operational nature, e.g., miles-in-trail, vectoring, etc. Indeed, these options are not considered in the tactical phase. In view of these modelling assumptions, the model captures delays due to flying longer paths and ground holding, but does not capture airborne holding delays due to air traffic control actions. 

The problem we are addressing is bi-objective in nature. There is a trade-off between the (total) preference score and the (total) delay assigned, as evidenced in \S \ref{sec::res}. In this work, our aim is to compute the efficient solution of the problem that either minimizes the total delay or allows delays to a given  threshold set by the relevant stakeholder. 

Finally, to improve the degree of ``realism'' of the model, we also provide an accurate representation of airspace capacities as detailed in the sequel.

\subsubsection{Airspace capacities}

All elements of the air traffic system, i.e., airports and en-route sectors, have limited capacity, which is an indirect measure of the maximum workload of an air traffic controller. Capacity is the maximum number of operations that can be safely executed within the element in a given interval of time. For the sake of accuracy, different definitions of the en-route sector capacity are used in Europe and in the US. In Europe, the en-route sector capacity limits the number of flights entering the sector in a given time interval (throughput), whilst in the US it limits the number of flights simultaneously in the sector (instant capacity). In what follows and without loss of generality, we adopt the European definition as we are going to apply our approach to instances of the European air traffic system. 
Capacities often change over time as they are affected by several factors including weather forecast and airspace configuration. Moreover, en-route sectors may have several values of the capacity at any point in time, each referring to a specific time interval. 
\begin{figure}
    \centering
    \includegraphics[width=0.9\textwidth]{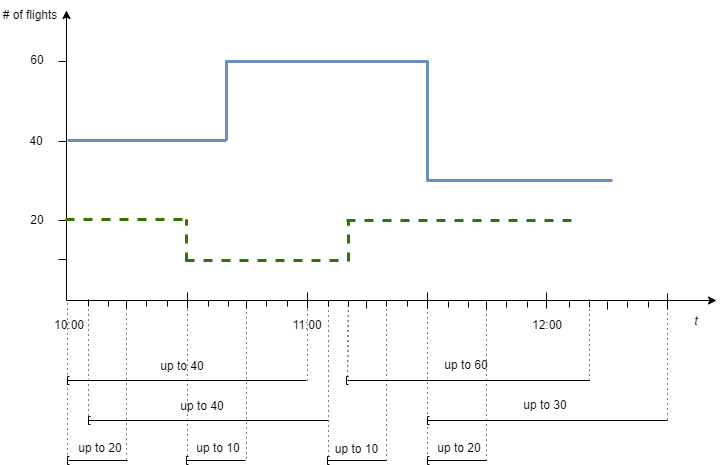}
    \caption{Sector capacity profiles and restrictions}
    \label{fig:Cap1}
\end{figure}
For instance, considering the illustrative example depicted in Figure \ref{fig:Cap1}, the sector has two values of the capacity. One for a one-hour time interval (solid line), and the second for a 15-minute time interval (dashed line). At 10 am, the one-hour capacity is 40, meaning that in the one-hour interval starting at 10 am, at most forty flights are allowed entering the sector. The same is for any one-hour time interval starting between 10:00 and 10:40 (excluded), when the hourly capacity takes value of 60. 
At 10 am and until 10:30 (excluded), the 15-minute capacity is 20, meaning that in any 15-minute time interval starting between 10 and 10:30 am (excluded), at most 20 flights are allowed entering the sector. This second value of the capacity, on a shorter interval of time, has the functionality of smoothing the sector air traffic demand out. Overall, these  capacity constraints provide more flexibility to the system  with respect to the usual modelling approach adopted in the literature, which normally limits the number of operations at any (discretized) period of time. This more realistic capacity representation allows for periods of (moderate) peak of the air traffic demand as long as these periods are preceded and/or followed by periods of low air traffic demand. As a by-product, this representation of the capacity may result in less conservative ATFM solutions, i.e., in a reduced use of control actions.

\subsubsection{The integer programming formulation}

To describe the trajectory-based model for ATFM, we introduce the following notation.

\subsubsection*{Sets}
\begin{description}
	\item[-] $T$: set of discrete time periods $ [t, t+1) $ dividing the time horizon; 
	\item[-] $S$: set of air traffic system elements, i.e., airports and en-route sectors; 
	\item[-] $F$: set of flights; 
	\item[-] $P_f$: set of 4D trajectories for flight $f \in F$ (computed by Algorithm \ref{alg:extrTraj});
	\item[-] $P_f(s)$: subset of 4D trajectories of $P_f$ crossing element $s \in S$; 
	\item[-] $I_{t,h} = \{t, t+1, \ldots, t+h-1\}$: subset of $h(\geq1)$ contiguous time periods starting at time period $t \in T$; 
	\item[-] $B_f$: set of ground delays that can be assigned to flight  $f \in F$; 
	\item[-] $H_{s,t}$: set of capacity limits for element $s \in S$ at time period $t \in T$. $h \in H_{s,t}$ is the duration of the capacity limit in number of discrete time periods.   
\end{description}

\subsubsection*{Parameters}
\begin{description}
\item[-] $G_p^f$: preference of flight $f \in F$ for trajectory $p \in P_f$ (computed by Algorithm \ref{alg:fliPref});
\item[-] $D_{pd}^f$: total delay suffered by flight $f \in F$ flying trajectory $p \in P_f$ with ground delay $d \in B_f$. The value of this parameter is  computed by 
	$ D_{pd}^f = \max \left(0 \; , \; {STD}^f + d + E_p - {STA}^f\right),
	$
	where  $STD^f$ ($STA^f$) is the scheduled time of departure (arrival) of flight $f \in F$; and  $E_p$ is the duration (in time periods) of trajectory $p \in P_f$.
\item[-] $R$: delay budget;
\item[-] $C_s^h(t)$: capacity of sector $s \in S$ at time period $t \in T$ for a {\it h}-period time interval;
\item[-] $\tau^f_{s,p}$: time period of entrance of flight $f \in F$ in sector $s \in S$ along trajectory $p \in P_f$ if leaving on time.
\end{description}

\subsubsection*{Decision variables}

$$ y_{pd}^f = \left\{ \begin {array}{clcr}
          1, & \mbox{if flight $f \in F$ flies trajectory $p \in P_f$ with ground delay $d \in B_f$,}\\
          0, & \mbox{otherwise.}
          \end {array}
\right.$$

\subsubsection*{Formulation}

\begin{alignat}{4}
\text{\emph{IP-pref}:} \qquad && \max \;  \sum_{f \in F} \sum_{p \in P_f} \sum_{d \in B_f} G_{p}^{f} \cdot y_{pd}^{f} \qquad & \\
&&\sum_{p \in P_f} \sum_{d \in B_f} y_{pd}^{f} = 1 \qquad  &\forall \; f \in F \label{cons::assign}\\
&&\sum_{f \in F, p \in P_f(s)} \; \sum_{d \in B_f : (\tau^f_{s,p} + d) \in I_{t,h} } y_{pd}^{f} \le C^h_s(t) \qquad  &\forall \; s \in S , t \in T, h \in H_{s,t} 
\label{cons::cap}\\
&&\sum_{f \in F} \sum_{p \in P_f} \sum_{d \in B_f} D^f_{pd} \cdot y_{pd}^{f} \le R \qquad &&\label{cons::dbudget}\\
&&y_{pd}^{f} \in \{0,1\} \qquad &\forall \; f \in F, p \in P_f, d \in B_f.
\label{cons::int} \end{alignat}

Constraints \eqref{cons::assign} impose that each flight flies one trajectory with possibly a ground delay assigned in departure. Constraints \eqref{cons::cap} formulate the capacity limits used in practise. Notice that the special case $H_{s,t} = \{1\}$ implements the US definition of instant capacity.
Finally, constraint \eqref{cons::dbudget} imposes an upper bound on the total delay assigned. If we are interested in computing the minimum delay efficient solution, then the delay budget will be set to the value computed by solving the following problem 
\begin{alignat}{4}
\text{\emph{IP-delay}:} \qquad \qquad \qquad && \min\; \left\{\sum_{f \in F} \sum_{p \in P_f} \sum_{d \in B_f} D^f_{pd}\cdot y_{pd}^{f} \hspace{2mm} : \hspace{2mm} \eqref{cons::assign}, \eqref{cons::cap}, \eqref{cons::int} \right\} & .
\label{Prod::delay}
\end{alignat}

\section{Solution techniques}
\label{sec::soltech}

Our ambition is to solve realistic daily instances of the ATFM problem. However, the number of trajectories extracted by Algorithm \ref{alg:extrTraj} is extremely large, order of half-dozen millions.  If these trajectories were inputed into the model, it would lead to very large scale instances of the problem with a number of variables largely exceeding  one hundred millions, which is not practicable. To address this issue, we first propose to reduce the number of trajectories feeding the mathematical model, see \S \ref{sec:sizeredclu}.
Second, we design a customized optimization approach to solve large scale instances of the problem. Indeed, even after reducing the number of trajectories, the daily ATFM instances are still prohibited for any solver and standard computer: a preliminary computational experience showed that a state-of-the-art solver is not able to compute the optimal solution for any of our benchmark instances due to computer’s memory limits that halted the solution process. Moreover, in several cases, the solver is not even able to compute a feasible solution after hours of computation, thus highlighting the need of a customized optimization technique. To overcome this computational challenge, in \S \ref{sec::dcg} we present a delayed column insertion approach 
to solve the linear relaxation. In addition to providing a dual bound of optimal integer solutions, this approach also provides a subset of variables to restrict the formulation to a convenient size.

\subsection{Reducing the model size by clustering techniques} 
\label{sec:sizeredclu}

In the mathematical formulation of \S \ref{sec::formulation}, trajectories are modelled in the {\it time-sector} space by binary vectors $A(p) \in \{0,1\}^{|S|\cdot |T|}$. The component of vector $A(p)$ corresponding to pair $(s, t) \in S \times T$ is equal to one if trajectory {\it p} intersects sector $s$ after $t$ time periods from departure, 0 otherwise. Due to time and 3D space discretization, some trajectories with similar but different 4D geometry may have the same representation in the {\it time-sector} space of the model. Therefore, from the optimization point of view, only one of these trajectories can be included in the model without affecting the set of feasible solutions and the objective function value. Indeed, these trajectories -- very likely -- belong to the same preference cluster, so they also have the same preference score.
In case of different preference scores of two trajectories, case that we did not observe in our computational experience, the trajectory with the highest score should be included in the model, as this one would be assigned to the flight. These observations justify the following remark.
\begin{remark}
The reduced model of {\it IP-pref} ({\it IP-delay})  -- obtained after discarding  trajectories that are equivalent in the {\it time-sector} space -- is equivalent to the original one. 
\end{remark}

To further reduce the number of variables, we select one representative trajectory for each group of similar trajectories in the time-sector representation. To this end, we apply a  $k$-mean clustering to the set of trajectories (in the time-sector representation) with a fairly large value of $k$. The shortest (in duration) trajectory of each cluster is chosen as a representative. In case of multiple trajectories with the same minimum duration, we select the closest one to the cluster center. Since the value of $k$ is fairly larger than the number of preference clusters output by Algorithm \ref{alg:extrTraj}, we observed that the trajectories of the same $k$-mean cluster have equal preference score. 

From now on,  \emph{IP-pref} and \emph{IP-delay} refer to the corresponding reduced-size models.

\subsection{A dynamic column insertion approach}
\label{sec::dcg}

To solve the problem, we use the optimization approach sketched in Algorithm \ref{alg:cgh}. The algorithm is composed of three main steps: {\it i}) computation of an initial solution of {\it IP-pref} or {\it IP-delay} by a fast greedy procedure; {\it ii}) solution of the linear relaxation by delayed column insertion; {\it iii}) calculation of the integer optimal solution of the formulation restricted to the subset of variables of the last restricted linear problem solved in step {\it ii}. 
\begin{algorithm}
\caption{\textsc{Dynamic column insertion for Trajectory Based ATFM} }
\label{alg:cgh}
\footnotesize

{\bf Input:} model IP (\emph{IP-pref} or \emph{IP-delay}), parameters $\alpha$ (aging) and $\beta$ (new columns).

{\bf Output:} a feasible solution to IP and related optimality gap.
\begin{algorithmic}[1]
\State Compute an initial solution 
to IP
\State \label{a:cgh:cgstart}Define model RLP as the linear relaxation of IP restricted to variables taking value 1 in the initial solution 
\Repeat 
    \State solve RLP 
    \State \label{a:cgh:redc} compute the reduced costs of the IP variables 
    \State \label{a:cgh:newcol} add to RLP the $\beta$ variables with smallest negative reduced cost
    \State \label{a:cgh:aging} remove from RLP variables that did not take part in the optimal basis during the last $\alpha$ iterations
\Until{no negative reduced-cost variable exists}\label{a:cgh:cgend}
\State solve the final RLP to integrality
\end{algorithmic}
\end{algorithm}

\subsubsection{Initial constructive heuristic.} \label{sec:initgreedy}

An initial solution of {\it IP-pref} or {\it IP-delay} is computed by a constructive heuristic that iteratively assigns a trajectory and, possibly, a ground delay to one flight at a time. The order with which the heuristic considers the flights is based on the scheduled time of departure, which may be considered as a proxy of the FPFS policy adopted by CASA. 
In case of delay minimization, to each flight $f$, we assign the $(p,d) \in P_f \times B_f$ pair that satisfies capacity constraints (\ref{cons::cap}) and corresponds to the highest preference score and lowest delay, in lexicographic order. This priority is intended to privilege most preferred trajectories in the initial solution before starting delay minimization. In case of tie, the pair that leaves the larger average per-cent residual sector capacity is selected. For preference maximization, the first two selection criteria -- i.e., preference and delay -- are swapped, in view of the total delay budget. Due to the imposed upper limit on the maximum delay that can be assigned to any flight, the greedy procedure may not be able to assign a $(p,d)$ pair to all flights without violating some capacity constraints. Therefore, one or more flights may remain unassigned and have to be fixed by the following steps of Algorithm \ref{alg:cgh}.

\subsubsection{Solving the linear relaxation by delayed column insertion.}

To solve the linear relaxation of the problem we use a delayed column insertion algorithm. The procedure iteratively solves a restricted linear program (RLP) with a subset of all the possible variables $y^f_{pd}$ (i.e., flight-trajectory-delay columns). The RLP is initialized with the subset of variables $y$ selected (i.e., set to one) by the initial greedy procedure and with possibly dummy variables, which are used to restore RLP feasibility in case of unassigned flights. 
At each iteration, the dual information of the RLP is used to search for negative reduced cost variables to be conveniently added (step 5). However, we impose a limit on the number of variables (parameter $\beta$) entering the RLP at each iteration (step \ref{a:cgh:newcol}). The reason for introducing such limit is because there is a trade-off between the number of iterations (convergence speed of the algorithm) and the computational effort required at each iteration due to the increasing size of the RLP. To restrain the rapid growth of RLPs, we discard variables of the RLP that are not likely to take positive value in the next optimal solution, since they have been out of basis in the last $\alpha$ -- aging parameter -- iterations (step \ref{a:cgh:aging}). By tuning these parameters, we guarantee a fast execution of the algorithm.

The column insertion algorithm stops as soon as no negative reduced cost variable exists, meaning that the solution to the current RLP is also optimal for the linear relaxation of the original non-restricted model, either {\it IP-pref} or {\it IP-delay}.

\subsubsection{Computing the integer solution.}

To compute a feasible integer solution, we solve the RLP with integrality constraints on decision variables back in place.  
The restricted integer linear program (RILP) so obtained is viable for standard solvers. Indeed, by properly calibrating parameters $\alpha$ and $\beta$ of the delayed column insertion algorithm, we obtain final RLPs of practicable size. This method can be regarded as a rounding scheme with the guarantee of providing the best solution within the space of ``rounded'' solutions. 
It is important to observe that random rounding schemes similar, for instance, to the one proposed in \cite{hamsa2014}, are not very effective for the specific problem herein considered. Indeed, they rarely provide a feasible solution, due to the lack of cancellation decisions in our model. 

We remark that the quality of the integer solution generated by the final step of Algorithm \ref{alg:cgh} can be assessed by comparing it to the bound provided by the optimal solution of the linear relaxation of {\it IP-pref} or {\it IP-delay}, as computed by the column insertion steps \ref{a:cgh:cgstart}-\ref{a:cgh:cgend}.

\section{Computational Experience}
\label{sec::CompExp}

In this Section, we present the computational experience with the data-driven approach of Sections \ref{sec::method} and \S \ref{sec::soltech} on a set of real numerical examples drawn from operational data sets. We used the \ECTRL~Demand Data Repository (DDR2) to extract data and feed both the predictive and the prescriptive analytics' components of our approach. DDR2 contains all relevant information for the description of the European Civil Aviation Conference (ECAC) area's air traffic system as well as historical data of the air traffic demand in the region.  The ECAC area's air traffic system is the second largest in the world per number of flight operations annually. It extends over Caucasia, Turkey and all the European countries with the exception of  Belarus and Russia. More specifically, DDR2 provides the information for the description of the 650 (on average) en-route sectors with the related activation times and capacities, functional air blocks (aggregations of en-route sectors), military areas and airports. 
It also provides data for flight identification (call sign, arrival/departure airports and times, aircraft type etc.) and related flight plans, i.e.,  4D trajectories for each flight.

The examples herein considered include air traffic movements (either departures or arrivals) from 916 airports across all the 44 member states of ECAC.

\subsection{ATFM instances}
\label{ssec::instances}

We consider ten instances of the ATFM problem each representing a whole day of operations in the ECAC area. The selected days correspond to the ten busiest days, i.e., with largest number of flights, during the summer of 2016. To avoid any air traffic demand pattern and/or bias, we selected at least one instance for each day of the week.
All the instances include more than 30000 flights, with the only exception of Saturday August 28, 2016 that had ``only" 28508 flights. Three of the instances (Fridays) exceed the 32,000 flights. 
Given the size of the instances, with thousands of origin-destination pairs, and the ensuing need of extremely large data set for the accurate prediction of preferences and trajectories, we partitioned all the flights into the three following categories: 
\begin{itemize}
    \item category A, which includes all the flights whose  airports of departure and arrival are both in the top twenty busiest airport sites -- i.e., with the largest number of movements -- of the ECAC area in the considered day. The involved airport sites are: London, Paris, Amsterdam, Istanbul, Frankfurt, Madrid, Barcelona, Palma de Mallorca, Rome, Berlin, Milan, Copenhagen, Munich, Dussendorlf, Zurich, Brussels, Stockholm, Athens, Dublin, Manchester, Vienna.  
    \item category B, which includes all the flights whose trajectories are entirely within the ECAC airspace and are not included in category A;
    \item category C, which includes all the remaining flights, i.e., flights directed to or coming from an airport outside the ECAC area. This category also contains flights that only fly over the ECAC area.
\end{itemize} 

For each flight, Algorithm \ref{alg:extrTraj} and the procedure described in \S\ref{sec:sizeredclu} provide a set of alternative trajectories to fly from airport of departure to destination. We recall that, because the number of feasible trajectories drawn from the DDR2 database for a given flight $f$ is extremely large (often order of thousands), the maximum number of trajectory options (parameter $k$ of the clustering procedure) is limited to: order of dozens (depending on the original size of $P_f$), for flights in category A; four, for flights in category B; one, for flights in category C.
This classification of the flights, though dictated by practical considerations on the solution approach herein proposed, is also motivated by operational considerations. Indeed, many of the flights in category B are short haul flights, feeding hubs in hub-and-spoke operations. For this class of flights, not many alternative trajectories are recorded.
\begin{table}
  \centering
  \caption{Instances}
  \begin{tabular}{l|c c c |c c c| c} \hline \hline
    day & $|F|$ & \multicolumn{1}{c}{$|F_A|$}& \multicolumn{1}{c|}{$|F_B|$}& \multicolumn{1}{c }{$|P|$} & \multicolumn{1}{c}{$P_A$} & \multicolumn{1}{c|}{$P_B$ }& \multicolumn{1}{c}{\#vars}\\
    \hline
Fri 08/07/16 & 31825 & 3823 & 21227 & 148,897 & 17.7 & 3.8 & 4,379,775 \\
Fri 26/08/16 & 32007 & 3540 & 21595 & 143,466 & 17.2 & 3.8 & 4,158,875 \\
Sat 27/08/16 & 28508 & 3113 & 18490 & 125,515 & 17.4 & 3.9 & 3,726,950   \\
Sun 28/08/16 & 30090 & 3347 & 19605 & 132,416 & 17.0 & 3.9 & 3,921,525 \\
Mon 29/08/16 & 31287 & 3692 & 20807 & 145,019 & 17.6 & 3.9 & 4,183,425 \\
Tue 30/08/16 & 30990 & 3692 & 20491 & 143,887 & 17.6 & 3.8 & 4,170,200 \\
Wed 31/08/16 & 31282 & 3760 & 20701 & 145,954 & 17.7 & 3.8 & 4,300,300 \\
Thu 01/09/16 & 31489 & 3871 & 20704 & 147,275 & 17.6 & 3.8 & 4,358,825 \\
Fri 02/09/16 & 32128 & 3821 & 21503 & 149,795 & 17.6 & 3.8 & 4,325,025 \\
Fri 09/09/16 & 32053 & 3932 & 21489 & 151,384 & 17.6 & 3.8 & 4,391,375 \\ \hline \hline
  \end{tabular}
  \label{tab:20cinstances}
\end{table}
Table \ref{tab:20cinstances} summarizes the main features of each instance. In the first four column, we report the date of the instance, the total number of flights included in the instance ($|F|$), the number of flights in category A ($|F_A|$) and in category B ($|F_B|$) respectively. Columns 5, 6 and 7 display the total number of alternative trajectories ($|P|$) and the average number of alternative trajectories per flight in category A ($P_A$) and category B ($P_B$) respectively. Finally, the last column displays the number of binary decision variables of the (reduced-size) mathematical model, which also depends on the maximum allowed ground delay that we set to 120 minutes for all the flights.

In the spirit of providing accurate solutions of the ATFM problem, the time period of discretization is a crucial element. The smaller is the discretization time period, the higher is the accuracy of solutions. However, this comes at higher computational expenses. In this work, we used a 5-minute discretization period. A finer discretization period would provide a level of detail that is usually captured in the tactical/operational phase and not in the planning one. Indeed, very accurate plans can be jeopardized by the inherent uncertainty affecting air traffic operations.

\subsection{Results}
\label{sec::res}

We implemented the models and the proposed solution techniques in C++ using Cplex C-API.  We used Cplex 12.9 as optimization engine and run the tests on a workstation equipped with an Intel Xeon E-2176G processor with 6 cores at 3.7 GHz, and 16 GB RAM.

Concerning the solution of retricted integer linear problem (RILP), we observe that solvers' routines for cut generation may be ineffective, thus making more difficult to close the gap without relaying on extensive branching. Therefore, we set a 1\% optimality gap and a one-hour time limit to solve the RILP.

To evidence the flexibility of the proposed approach, we here report the computational performance to solve both the minimization of total delays (model \emph{IP-delay}) and the maximization of total preference scores (\emph{IP-pref}). Indeed, both the problems have to be solved to compute any efficient solution of the problem, including the one  minimizing the total delays, unless information on the magnitude of the two objectives' value is available. For each instance of the problem, we give statistics to compute the optimal solution of both the linear relaxation (LR) -- with the delayed column insertion algorithm -- and the RILP. For the delayed column insertion algorithm, we report the computational time ({\it Time}) in seconds, the optimal value of the linear relaxation problem ($z^{LR}$) and  the number of iterations ({\it Iter.s}). For the RILP, we report the number of variables ({\it \# var.s}) of the mathematical program, the solver's computational time ({\it Time}) in seconds, the value of the  optimal solution ($z^{\ast}$), and the percentage gap ({\it GAP}) between $z^{\ast}$ and the linear relaxation lower bound ($z^{LR}$). Finally,  we also report  the total time ({\it Total Time}) in minutes to solve the instance. 


\subsubsection*{Minimization of delays.}
Table~\ref{tab:delay} summarizes the results on \emph{IP-delay}. 
\begin{table}
\centering
  \caption{Computational performances for total delays' minimization.}
  \begin{tabular}{c |c c c | c c c  c | c } \hline \hline
 \multicolumn{1}{c}{}    & \multicolumn{3}{c}{LR} & \multicolumn{4}{c}{RILP} & \multicolumn{1}{c}{}   \\ \hline 
Instance	&	Time	&	$z^{LR}$	&	Iter.s	&	\# var.s	&	Time	&	$z^\ast$	&	GAP	&	Total Time		\\
	&	(secs)	&		&		&		&	(secs)	&		&	\%	&	(mins)		\\	\hline
08/07/2016	&	785	&	5771.1	&	51	&	44723	&	7	&	5785	&	0.24	&	13.20	\\	
26/08/2016	&	615	&	5926.1	&	48	&	44375	&	11	&	5929	&	0.05	&	10.43	\\	
27/08/2016	&	2033	&	7485.4	&	79	&	41643	&	350	&	7558	&	0.97	&	39.72	\\	
28/08/2016	&	1148	&	6293.6	&	51	&	42153	&	66	&	6315	&	0.34	&	20.23	\\	
29/08/2016	&	685	&	6293.0	&	47	&	44446	&	13	&	6309	&	0.25	&	11.63	\\	
30/08/2016	&	642	&	5125.8	&	43	&	45024	&	25	&	5137	&	0.22	&	11.12	\\	
31/08/2016	&	549	&	6045.0	&	49	&	43788	&	16	&	6075	&	0.50	&	9.42	\\	
01/09/2016	&	623	&	5843.2	&	52	&	45031	&	12	&	5847	&	0.07	&	10.58	\\	
02/09/2016	&	403	&	5883.5	&	47	&	43824	&	4	&	5891	&	0.13	&	6.78	\\	
09/09/2016	&	499	&	5678.2	&	47	&	44970	&	12	&	5727	&	0.86	&	8.52	\\	\hline
AVERAGE	&	798	&		&		&		&	52	&		&	0.36	&	14.16	\\
 \hline \hline
\end{tabular}
  \label{tab:delay}
\end{table}
The computational time to solve the linear relaxation is rather short. The  average ({\it resp.} median)  is 798.5 ({\it resp.} 632.5) seconds and only few dozen iterations are needed to converge to the optimal solution. The resolution of the RILP is also very fast. The average ({\it resp.} median) of the computational times is 51.6  ({\it resp.} 12.5) seconds, with only two instances --namely August $27^{th}$ and August $28^{th}$ -- exceeding one minute. The optimality gap is also very small with an average value of 0.36\%. It certifies that we are able to compute a solution that is very close to the optimal solution of {\it IP-delay} for all the instances, with a gap that is consistently lower than 1\%. 
But even more important, we are able to compute the optimal solution for all the benchmark instances in at most 39.7 minutes, which is acceptable from the practical point of view.  The average of the total computational times is only 14.2 minutes, despite the relatively long computational time of the August $27^{th}$ instance.

\subsubsection*{Maximization of the total preference scores.}
Table~\ref{tab::preferences} summarizes the results on model \emph{IP-pref} to compute the efficient solution for a total delay budget -- 
constraint \eqref{cons::dbudget} of the model's formulation --  set to 110\% of the minimum delay. Also for this problem, we are able to compute optimal solutions for all the benchmark instances. 
\begin{table}
\centering
  \caption{Computational performances for the maximization of the total preference scores.}
  \begin{tabular}{c |c c c | c c c c | c } \hline \hline
 \multicolumn{1}{c}{}    & \multicolumn{3}{c}{LR} & \multicolumn{3}{c}{RILP} &  \multicolumn{1}{c}{} \\ \hline 
Instance	&	Time	&	$z^{LR}$	&	Iter.s	&	\# var.s	&	Time	& $z^{\ast}$	&	GAP	&	Total Time			\\
	&	(secs)	&		&		&		&	(secs)	&		&	\%	&	(mins)		\\	\hline
08/07/2016	&	2052	&	24458.3	&	71	&	83710	&	48	&	24437.0	&	0.09	&	35.0	\\
26/08/2016	&	1104	&	24524.1	&	66	&	83845	&	20	&	24393.5	&	0.53	&	18.7	\\
27/08/2016	&	9472	&	21341.8	&	116	&	70550	&	75	&	21135.9	&	0.96	&	159.1	\\
28/08/2016	&	3731	&	22795.9	&	65	&	75968	&	31	&	22608.1	&	0.82	&	62.7	\\
29/08/2016	&	1406	&	23926.8	&	79	&	85860	&	48	&	23916.1	&	0.04	&	24.2	\\
30/08/2016	&	1442	&	23619.9	&	65	&	83845	&	54	&	23584.9	&	0.15	&	24.9	\\
31/08/2016	&	1157	&	23989.9	&	75	&	83474	&	22	&	23827.9	&	0.68	&	19.7	\\
01/09/2016	&	1064	&	24294.9	&	67	&	85273	&	21	&	24153.7	&	0.58	&	18.1	\\
02/09/2016	&	909	&	24615.2	&	75	&	87511	&	23	&	24478.5	&	0.56	&	15.5	\\
09/09/2016	&	866	&	24310.1	&	74	&	87315	&	31	&	24294.1	&	0.07	&	15.0	\\ \hline
AVERAGE	&	2320	&		&		&		&	37	&		&	0.45	&	39.29	\\
 \hline \hline
  \end{tabular}
  \label{tab::preferences}
\end{table}
The optimality gap of the RILP solution is 0.45\% on average, and never exceeds 1\%. However, we observe a moderate deterioration of the computational times with respect to the delay minimization experience. The average ({\it resp.} median) of the total computational times is 39.3 ({\it resp.} 21.9) minutes. The longer computational times have to be ascribed to the resolution of the linear relaxation. In fact, solving the RILP is consistently fast, with only one instance slightly exceeding one minute of computation, i.e., 75 seconds. On the other hand, the running time of the delayed column insertion is relatively long. The average ({\it resp.} median) value over the set of benchmark instances is 2320 ({\it resp.} 1595) seconds. The longer computational times are due to a larger number of iterations as well as to longer  time per iteration. 
The RLP grows faster than its corresponding problem in the case of delays minimization thus requiring a higher computational effort and time to be solved. However, it is also important to highlight that computational times are - by and large - consistent with the potential use of the model in practice. The are only two instances, i.e., August $27^{th}$ and August $28^{th}$, which require more than one hour of computation to be solved. As mentioned above, almost all of this time is used to solve the linear relaxation of the problem. The improvement of the RLP's objective function between two consecutive iterations of the delayed column insertion algorithm tends to vanish, thus requiring a relatively large number of iterations to converge to the optimal solution.
This behaviour justifies the use of time limits on the delayed column insertion algorithm without compromising much on the quality of the final solution.    
\begin{table}
\begin{small}
\centering
\caption{Integer feasible solutions with varying time limits for the August $27^{th}$ and $28^{th}$ instances.}
\begin{tabular}{c |c c | c c c || c c | c c c} \hline \hline
\multicolumn{1}{c}{}& \multicolumn{5}{c}{August $27^{th}$} & \multicolumn{5}{c}{August $28^{th}$}\\ \hline
T.Limit	&	$z^{LR}$	&	$GAP^{LP}$	&	$z^{IP}$	&	Time	&	$GAP^{IP}$	&	$z^{LR}$	&	$GAP^{LP}$	&	$z^{IP}$	&	Time	&	$GAP^{IP}$	\\ \hline
2000	&	21022.3	&	1.50	&	20948.9	&	37.5	&	1.84	&	22719.5	&	0.34	&	22499.1	&	36.9	&	1.30	\\
2500	&	21104.4	&	1.11	&	21031.3	&	44.8	&	1.46	&	22750.2	&	0.20	&	22542.2	&	42.6	&	1.11	\\
3000	&	21160.3	&	0.85	&	20980.3	&	53.0	&	1.69	&	22781.3	&	0.06	&	22568.7	&	51.4	&	1.00	\\
3400	&	21181.9	&	0.75	&	20971.5	&	59.4	&	1.74	&	22793.4	&	0.01	&	22588.5	&	58.6	&	0.91	\\
No Limit	&	21341.8	&	0.00	&	21135.9	&	159.1	&	0.96	&	22795.9	&	0.00	&	22608.1	&	62.7	&	0.82	\\
\hline \hline
  \end{tabular}
  \label{tab::trendLP}
\end{small}
\end{table}
In Table \ref{tab::trendLP}, we summarize the values of the RLP's solution ($z^{LR}$) and the corresponding integer solution ($z^{IP}$) using different time limits. The table also shows the accuracy of the RLP's solution ($GAP^{LP}$) and the optimality gap of the solution of the corresponding RILP ($GAP^{IP}$), together with the total (delayed column insertion and integer linear programming solver) running time ({\it Time}) in minutes. For the August $27^{th}$ instance, with a time limit of 2000 seconds, we compute an integer feasible solution whose optimality gap is only 1.84\%, which is reasonable for the use in practise. Moreover,  this solution is not even 1\% worse than the best integer solution computed  with no time limit. Even better results are obtained for the August $28^{th}$ instance, where the optimality gap after 2000 seconds is 1.30\%, which is not even 0.5\% larger than the optimality gap without any time limit.

The following remark summarizes the key contribution of the  results described above. 

\begin{remark}
\label{rem::opt}
The proposed approach computes solutions of the ATFM problem that are optimal or very close to the optimum for all the instances we tested. Computational times are short and consistent with the potential use of the method in practice.
\end{remark}

\subsubsection*{Analysis of solutions and the delay-preference trade-off.}
\begin{figure}
     \centering
     \begin{subfigure}[]{0.49\textwidth}
         \centering
         \vspace{0.92cm}
         \includegraphics[width=\textwidth]{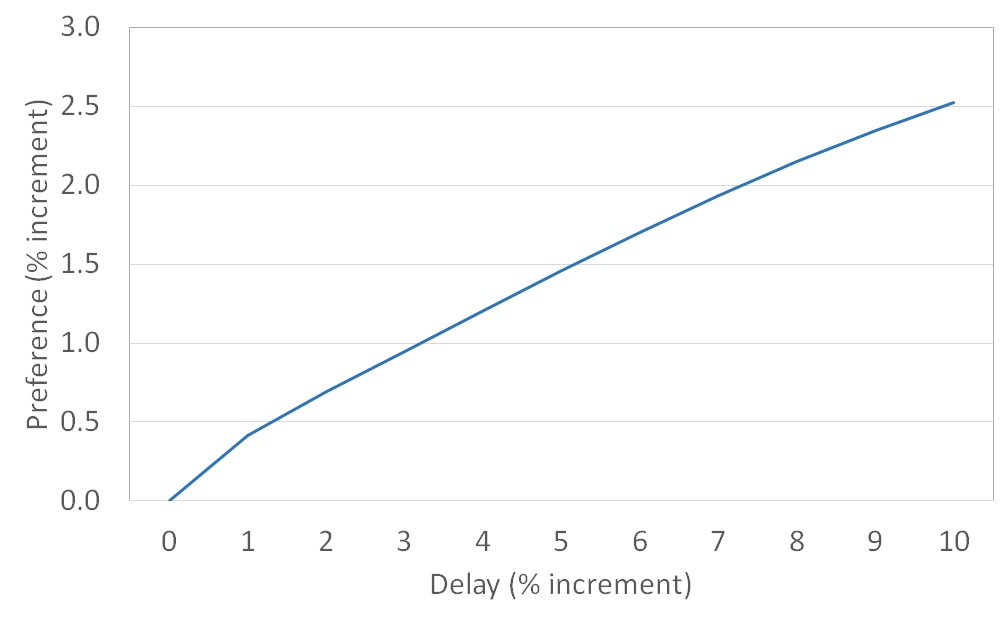}
         \caption{\footnotesize Pareto efficient frontier}
    \label{fig::example}
     \end{subfigure}
     \hfill
     \begin{subfigure}[]{0.49\textwidth}
         \centering
         \includegraphics[width=\textwidth]{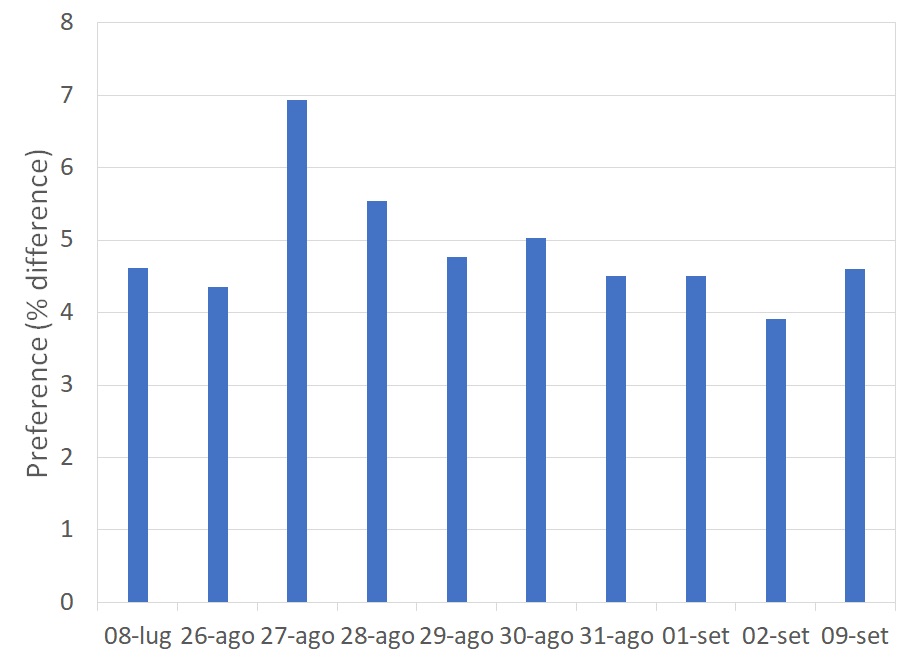}
         \caption{\footnotesize ``Inefficency'' of the min delay optimal solution.}
         \label{fig::innef}
     \end{subfigure}
   \caption{\it On the delays -- preference scores trade-off.}
    \label{fig:results_existing_methods}
\end{figure}

To conclude this last section, we highlight the trade-off between (total) delays and preference scores. Figure \ref{fig::example} displays an approximation of the Pareto frontier -- i.e., the set of efficient solutions -- of September 2$^{nd}$ instance, obtained considering the linear relaxation of the problem. We report delays on the abscissa and the preference scores on the ordinate. More in particular, each axis displays the percentage increment with respect to the efficient solution’s value that minimizes the total delay. 
In Figure \ref{fig::innef}, the bar chart represents the ``inefficiency'' value of the optimal solutions of problem {\it IP-delay}. The inefficiency value is the percentage difference of the total preference scores of the delay minimizing optimal solution and the corresponding efficient solution.
As highlighted by the bar chart, none of the delay minimizing optimal solutions is Pareto efficient. The average inefficiency value is 4.88\%,  with a peak of almost 7\% for the August $27^{th}$ instance. Indeed, a 5\% increase of the preference score has the practical implication that at least 1500 flights --  a very conservative estimate -- can fly a trajectory with a higher preference. Potentially, the inefficiency value can be much larger than the ones observed, as it is possible to compute delay minimizing solutions with an inefficiency value greater than 50\% for all the benchmark instances.
This demonstrates the importance of computing efficient solutions for the ATFM problem, as these solutions can be more likely to be accepted by the stakeholders, thus facilitating the achievement of a consensus within the decision-making process.

\section{Conclusions}
\label{sec::concl}

In this paper, we have presented a data-driven optimization approach for ATFM trajectory based operations that makes three significant contributions. First, it provides an accurate representation of air traffic system’s operations, by using both flown 4D trajectories as well as a refined modelling of airspace system’s capacity constraints. Second, it proposes a machine learning approach to assess airspace users’ preferences, thus overcoming some of the barriers that have been experienced in practise. The third and, from the practical viewpoint, most important contribution of the model is that it is computationally viable even for daily instances of the problem with up to 32,000 flights. 

The computational experiments suggest that ATFM problems of a size comparable to the entire network of en-route sectors and airports in the European air traffic system can be solved to near-optimality within reasonable -- for the applications context -- computation times. 
Indeed, running times of the order of 30 minutes are consistent with the time constants associated with the current decision cycles at the \ECTRL's Network Manager, the facility that coordinates ATFM for the entire European airspace. Moreover, the computational experiments 
%
show the importance of considering preferences, since significant differences can be observed in terms of overall preference score between equivalent solutions from the delay minimization standpoint. By computing Pareto efficient solutions, we expect to have a smaller number of iterations between the Network Manager and airspace users and a faster and fine-tuned resolution of demand-capacity imbalances.

Our long term ambition is to feed the proposed approach into industrial research activities and to enable a derivative version of this model to become one of the basic network management decision support tools of the future air traffic management system.

\bibliography{main_DLL}

\end{document}